\documentclass{article}
\begin{document}
%
%
%
%
%
\title{Imaginary eigenvalues
and complex eigenvectors explained by real geometry}
\author{Eckhard M.S. Hitzer\footnote{Department of Mechanical Engineering, Fukui University, 
3-9-1 Bunkyo, 910-0024 Fukui, Japan, e-mail: hitzer@mech.fukui-u.ac.jp}}
\date{11 July 2001}
\maketitle

\begin{abstract}
This paper first reviews how anti-symmetric matrices in two dimensions yield imaginary
eigenvalues and complex eigenvectors. It is shown how this carries on to rotations by
means of the Cayley transformation. Then
a real geometric interpretation is given to the eigenvalues and eigenvectors
by means of real geometric algebra.
The eigenvectors are seen to be \textit{two component eigenspinors} which can be further reduced to underlying
vector duplets. The eigenvalues are interpreted as rotation operators, which rotate the underlying vector
duplets. The second part of this paper extends and generalizes the treatment to three dimensions.
Finally the four-dimensional problem is stated.
\end{abstract}

\section{Introduction}

\begin{quotation}
\ldots for geometry, you know, is the gate of science, and the gate
is so low and small that one can only enter it as a little child.
(\textit{William K. Clifford}~\cite{Gu:inotr})
\end{quotation}
\begin{quotation}
 But the gate to life is narrow and the way that leads to it is hard,
and there are few people who find it. \ldots I assure you that unless you change
and become like children, you will never enter the Kingdom of heaven.
(\textit{Jesus Christ}~\cite{JC:bible})
\end{quotation}

This article arose from a linear algebra problem on anti-symmetric matrices for undergraduate engineering
students. I wrote it, looking for a real geometric understanding of the imaginary eigenvalues and complex 
eigenvectors. Being already familiar with geometric algebra\cite{Gu:inotr,DH:STC,DH:CtoG,DH:NF1} it was natural
to try to apply it in this context. I haven't come to terms with the four dimensional real interpretation,
but I think it worthwhile at the end, to at least state the (complex) problem. 

The first quotation stems from Clifford himself, who initially was a theologian and then became an atheist. 
But somehow his view of science was strongly colored by what Jesus taught as the Gospel about the 
Kingdom of God.\footnote{It is interesting to note that a parallel to this exists even in the Japanese tea ceremony: 
"\ldots enter the teahouse. The sliding door is only thirty six inches high. Thus all who enter must bow their heads and crouch. 
This door points to the reality that all are equal in tea, irrespective of status or social position."\cite{URL:teac}. 
The present form of the tea ceremony was established by Sen Rikyu in the 16th century. His wife is said to have been
a secret Christian (Kirishitan), some think even Sen Rikyu was. } 
To agree or disagree on what Clifford believed is a matter of faith and not of science. 
But I quite like his point, that geometry is like a gateway to a new understanding of 
science.

\section{Two real dimensions}
\subsection{Complex treatment}

Any anti-symmetric matrix in two real dimensions is proportional to
\begin{equation}
\label{eq:U2}
U=\left(
	\begin{array}{cc}
	0 & -1 \\
	1 & 0	
	\end{array}
	\right).
\end{equation}
The characteristic polynomial equation of the matrix $U$ is
\begin{equation}
| U-\lambda E | = \left|
                      \begin{array}{cc}
	              -\lambda & -1 \\
	              1 & -\lambda	
	              \end{array}
                      \right|
= \lambda^{2}+1
=0,
 \;\;\;\;\;
\mbox{i.e.}\;\;\;\;\;
\lambda^{2}=-1.
\label{eq:cpe2}
\label{eq:cpe2a}
\end{equation}
The classical way to solve this equation is to postulate an imaginary entity $j$
to be the root
$j=\sqrt{ -1}$. This leads to many interesting consequences, yet any real
geometric meaning of this imaginary quantity is left obscure.  
The two eigenvalues are therefore the imaginary unit $ j$ and $- j$ .
\begin{equation}
\lambda_{1}=j, \,\,\, \lambda_{2}=-j
\end{equation}
The corresponding complex eigenvectors $\mathbf{x}_1$ and $\mathbf{x}_2$ are
\begin{equation}
U \mathbf{x}_1= \lambda_1 \mathbf{x}_1 = j \mathbf{x}_1 \rightarrow
 \mathbf{x}_1=\left(
	\begin{array}{c}
	1 \\
	-j	
	\end{array}
	\right),
\;
U \mathbf{x}_2= \lambda_2 \mathbf{x}_2 = - j \mathbf{x}_2 \rightarrow
 \mathbf{x}_2=\left(
	\begin{array}{c}
	1 \\
	j	
	\end{array}
	\right).
\end{equation}

The Cayley transformation~\cite{Na:senkei} $C(-kU)$, with $k=(1-\cos \vartheta)/\sin \vartheta$
\begin{eqnarray}
C(-kU) & = & (E+(-kU))^{-1}{}(E-(-kU))
=  E-\frac{2}{1+k^2}(-kU-(-kU)^2) \nonumber \\
& = & \left(
	\begin{array}{cc}
	\cos \vartheta & -\sin \vartheta  \\
	\sin \vartheta & \cos \vartheta	
	\end{array}
	\right)
\label{eq:C2}
\end{eqnarray}
allows to describe two-dimensional rotations.

The third expression of equation (\ref{eq:C2}) shows that $ U$
and $C(-kU)$ must have the same eigenvectors $\mathbf{x}_1$ and $\mathbf{x}_2$.
The corresponding eigenvalues of $C(-kU)$ can now easily be calculated from (\ref{eq:C2}) as
\begin{equation}
\lambda_{c1}=1+\frac{2k\lambda_1(1+k\lambda_1)}{1+k^2}, \;\;\;\;\;
\lambda_{c2}=1+\frac{2k\lambda_2(1+k\lambda_2)}{1+k^2}.
\end{equation}
Inserting $\lambda_{1}=j$ and $ \lambda_{2}=-j$
we obtain the complex eigenvalues of the two-dimensional rotation $C(-kU)$ as
\begin{equation}
\lambda_{c1}=\cos \vartheta + j \sin \vartheta, \;\;\;\;\;
\lambda_{c2}=\cos \vartheta - j \sin \vartheta.
\end{equation}

We now face the question what the imaginary and complex eigenvalues and the complex
eigenvectors of $U$ and the rotation $C(-kU)$ mean in terms of purely real geometry.
In order to do this let us turn to the real geometric algebra $\mathrm{R}_2$
of a real two-dimensional vector space $\mathrm{R}^2$.~\cite{Gu:inotr,DH:CtoG,DH:NF1}

\subsection{Real explanation}

\label{sc:2rex}
Instead of postulating the imaginary unit $j$ we now solve the characteristic polynomial
equation (\ref{eq:cpe2}) using both orientations of the oriented unit area element
\textbf{i} of $\mathrm{R}_2$:
\begin{equation}
\lambda_1=\mathbf{i}, \;\;\; \lambda_2=-\mathbf{i}.
\end{equation}
The corresponding "eigenvectors" $\mathbf{x}_1$ and $\mathbf{x}_2$ will then be:
\begin{equation}
\mathbf{x}_1=\left(
	\begin{array}{c}
	1 \\
	-\mathbf{i}
	\end{array}
	\right), \;\;\;
\mathbf{x}_2=\left(
	\begin{array}{c}
	1 \\
	\mathbf{i}
	\end{array}
	\right).
\end{equation}

As before, the "eigenvectors" of the Cayley transformation $C(-kU)$ will be the same.
And the eigenvalues of $C(-kU)$ now become:
\begin{equation}
\lambda_{c1}=\cos \vartheta + \mathbf{i} \sin \vartheta, \;\;\;\;\;
\lambda_{c2}=\cos \vartheta - \mathbf{i} \sin \vartheta.
\end{equation}

We can now take the first step in our real explanation and identify
the two "eigenvectors" $\mathbf{x}_1$ and $\mathbf{x}_2$ as \textit{two-component spinors} with
the entries: $x_{11}=1, \; x_{12}= - \mathbf{i}$ and  $x_{21}=1, \; x_{22}= \mathbf{i}$.

Now we want to better understand what the real-oriented-unit-area-element eigenvalues
$\lambda_1$,  $\lambda_2$ as well as
$\lambda_{c1}$ and 
$\lambda_{c2}$ do when multiplied with the two-component
eigen-spinors $\mathbf{x}_1$ and $\mathbf{x}_2$. Every spinor can
be understood to be the geometric product of two vectors. We therefore choose an arbitrary,
but fixed reference vector \textbf{z} from the vector space $\mathrm{R}^2$.
For simplicity let us take \textbf{z} to be $\mathbf{z}=\mbox{\boldmath $\sigma$}_1$, assuming
\{$\mbox{\boldmath $\sigma$}_1,\mbox{\boldmath $\sigma$}_2$\} to be the orthonormal basis of $\mathrm{R}^2$.
We can then factorize the spinor components of the eigen-spinors $\mathbf{x}_1$ and $\mathbf{x}_2$ to:
\begin{equation}
x_{11}=x_{21}=1=\mbox{\boldmath $\sigma$}_1\mbox{\boldmath $\sigma$}_1,\;\;
x_{12}=-\mathbf{i}=\mbox{\boldmath $\sigma$}_2\mbox{\boldmath $\sigma$}_1,\;\;
x_{22}=\mathbf{i}=-\mbox{\boldmath $\sigma$}_2\mbox{\boldmath $\sigma$}_1
\end{equation}
Note that we always factored out $\mathbf{z}$ to the right. \label{pg:2rightf}
In two real dimensions it now seems natural to adopt the following interpretation: 
The eigen-spinor $\mathbf{x}_1$ corresponds (modulus the geometric 
multiplication from the right with $\mathbf{z}=\mbox{\boldmath $\sigma$}_1$) to the real vector pair 
\{$\mbox{\boldmath $\sigma$}_1,\mbox{\boldmath $\sigma$}_2$\}, 
whereas $\mathbf{x}_2$ corresponds to the real vector pair 
\{$\mbox{\boldmath $\sigma$}_1,-\mbox{\boldmath $\sigma$}_2$\}. 
Multiplication with $\lambda_1$ from the left as in
\begin{equation}
U\mathbf{x}_1=\lambda_1\mathbf{x}_1=\mathbf{i}\left(
	\begin{array}{c}
	x_{11} \\
	x_{12}
	\end{array}
	\right)
=\left(
	\begin{array}{c}
	\mathbf{i}x_{11} \\
	\mathbf{i}x_{12}
	\end{array}
	\right)
\end{equation}
results in
\begin{equation}
x_{11} \rightarrow \mathbf{i} x_{11}=(-\mbox{\boldmath $\sigma$}_2)\mbox{\boldmath $\sigma$}_1,\;\;\;
x_{12} \rightarrow \mathbf{i} x_{12}=\mbox{\boldmath $\sigma$}_1\mbox{\boldmath $\sigma$}_1
\end{equation}
That is the multiplication with $\lambda_1 =\mathbf{i}$ from the left transforms the
vector pair \{$\mbox{\boldmath $\sigma$}_1,\mbox{\boldmath $\sigma$}_2$\} 
to the new pair \{-$\mbox{\boldmath $\sigma$}_2,\mbox{\boldmath $\sigma$}_1$\}, 
which is a simple rotation by $-90$ degrees. 
Here the non-commutative nature of the geometric product is important. 

The analogous calculation for $\lambda_2\mathbf{x}_2=-\mathbf{i}\mathbf{x}_2$ shows that the pair 
\{$\mbox{\boldmath $\sigma$}_1,-\mbox{\boldmath $\sigma$}_2$\}, 
which corresponds to $\mathbf{x}_2$ is transformed to 
\{$\mbox{\boldmath $\sigma$}_2,-\mbox{\boldmath $\sigma$}_1$\}, i.e. it is rotated by +90 degree. 

I will now treat $C(-kU)\mathbf{x}_1=\lambda_{c1}\mathbf{x}_1$ 
and $C(-kU)\mathbf{x}_2=\lambda_{c2}\mathbf{x}_2$ in the same way.
\begin{equation} 
x_{11}\rightarrow \lambda_{c1} x_{11}
=(\cos \vartheta + \mathbf{i} \sin \vartheta)\mbox{\boldmath $\sigma$}_1^2
=\mbox{\boldmath $\sigma$}_1(\cos \vartheta - \mathbf{i} \sin \vartheta)\mbox{\boldmath $\sigma$}_1
=(\mbox{\boldmath $\sigma$}_1 R(-\vartheta))\mbox{\boldmath $\sigma$}_1,
\end{equation}
where $R(-\vartheta)$ is the rotation operator by $-\vartheta$. For the second component $x_{12}$
we have
\begin{equation} 
x_{12}\rightarrow \lambda_{c1} x_{12}
=(\cos \vartheta + \mathbf{i} \sin \vartheta)(\mbox{\boldmath $\sigma$}_2\mbox{\boldmath $\sigma$}_1)
=(\mbox{\boldmath $\sigma$}_2 R(-\vartheta))\mbox{\boldmath $\sigma$}_1.
\end{equation}
The action of $\lambda_{c1}$ on $\mathbf{x}_1$ means therefore a rotation of the corresponding vector pair
\{$\mbox{\boldmath $\sigma$}_1,\mbox{\boldmath $\sigma$}_2$\}  by $-\vartheta$.

The analogous calculations for $\lambda_{c2}\mathbf{x}_2$ show that $\lambda_{c1}$ rotates the vector pair 
\{$\mbox{\boldmath $\sigma$}_1,-\mbox{\boldmath $\sigma$}_2$\}, which corresponds to $\mathbf{x}_2$, into 
\{$\mbox{\boldmath $\sigma$}_1 R(\vartheta),-\mbox{\boldmath $\sigma$}_2 R(\vartheta)$\}. This corresponds to a
rotation of the vector pair by $+\vartheta$.

Summarizing the two-dimensional situation, we see that the complex eigenvectors $\mathbf{x}_1$ and $\mathbf{x}_2$ may 
rightfully be interpreted as two-component 
eigen-spinors with underlying vector pairs. The multiplication of these eigen-spinors with the 
unit-oriented-area-element eigenvalues $\lambda_1$ and $\lambda_2$ means a real rotation of the underlying vector 
pairs by -90 and +90 degrees, respectively. 
Whereas the multiplication with $\lambda_{c1}$ and $\lambda_{c2}$ means a real rotation of the underlying 
vector pairs by $-\vartheta$  and $+\vartheta$, respectively.

Now all imaginary eigenvalues and complex eigenvectors of anti-symmetric matrices in two real dimensions have
a real geometric interpretation. Let us examine next how this carries on to three dimensions.

\section{Three real dimensions}

\subsection{Complex treatment of three dimensions}

\label{sc:3comp}
Any anti-symmetric matrix in three real dimensions is proportional to a matrix of the form
\begin{equation}
U=\left(
\begin{array}{ccc}
0 & -c & b \\
c & 0 & -a \\
-b & a & 0
\end{array}
\right)
\end{equation}
with $a^2+b^2+c^2=1$. The characteristic polynomial equation of the matrix $U$ is
\begin{equation}
|U-\lambda E| = 
\left|
\begin{array}{ccc}
-\lambda & -c & b \\
c & -\lambda & -a \\
-b & a & -\lambda
\end{array}
\right|
= \lambda (\lambda^2+a^2+b^2+c^2)
=0.
\end{equation}
If we use the condition that $a^2+b^2+c^2=1$, this simplifies and breaks up into the two equations
\begin{equation}
\lambda^2_{1,2}=-1, \;\;\; 
\lambda_3 = 0.
\label{eq:cpe3}
\end{equation}

That means we have one eigenvalue $\lambda_3$ equal to zero and for the other two eigenvalues $\lambda_1$, $\lambda_2$ 
we have the same condition (\ref{eq:cpe2a}) as in the two-dimensional case for the matrix of equation (\ref{eq:U2}). 
It is therefore clear that in the conventional treatment one would again 
assign\footnote{Here an incompleteness of the conventional treatment becomes obvious. A priori there is no
reason to assume that the solutions to the characteristic polynomial equations in two and three dimensions 
(\ref{eq:cpe2a}) and (\ref{eq:cpe3})
must geometrically be the same.} 
$\lambda_1=j$ and $\lambda_2=-j$. 
The corresponding complex eigenvectors are:
\begin{equation}
\mathbf{x}_1=
\left( \begin{array}{c}
        1-a^2 \\
        -ab-jc \\
        -ac+jb
       \end{array}\right)
\doteq
\left( \begin{array}{c}
        -ab+jc \\
        1-b^2 \\
        -bc-ja
       \end{array}\right)
\doteq
\left( \begin{array}{c}
        -ac-jb \\
        -bc+ja \\
        1-c^2
       \end{array}\right), \;\;\;
\mathbf{x}_2= \mathrm{cc}(\mathbf{x}_1),
\label{eq:cevec3}
\end{equation}
where cc(.) stands for the usual complex conjugation, i.e. $ \mathrm{cc}(j)=-j$. The symbol $\doteq$ expresses 
that all three given forms are equivalent up to the multiplication with a scalar (complex) constant. 

The eigenvector that corresponds to $\lambda_3$ simply is:
\begin{equation}
\mathbf{x}_3=
\left( \begin{array}{c}
        a \\
        b \\
        c
       \end{array}\right).
\end{equation}
The fact that $\lambda_3=0$ simply means that the matrix $U$ projects out any component 
of a vector parallel to $\mathbf{x}_3$. 
$U$ maps the three-dimensional vector space therefore to a plane perpendicular to 
$\mathbf{x}_3$ containing the origin. 

The Cayley transformation~\cite{Na:senkei} $C(-kU)$ with $k=\frac{1-\cos \vartheta}{\sin \vartheta}$
now describes rotations in three dimensions:
\begin{eqnarray}
C(-kU)& =& (E+(-kU))^{-1}{}(E-(-kU)) \nonumber \\
& = & E-\frac{2}{1+k^2(a^2+b^2+c^2)}(-kU-(-kU)^2)= 
\label{eq:CayM3}
\end{eqnarray}
$$ \left(
	\begin{array}{ccc}
	1+(1-\cos \vartheta)(1-a^2) & -c\sin \vartheta+ab(1-\cos \vartheta) & b\sin \vartheta+ac(1-\cos \vartheta) \\
	c\sin \vartheta+ab(1-\cos \vartheta) & 1+(1-\cos \vartheta)(1-b^2) & -a\sin \vartheta+bc(1-\cos \vartheta) \\
        -b\sin \vartheta+ac(1-\cos \vartheta) & a\sin \vartheta+bc(1-\cos \vartheta) & 1+(1-\cos \vartheta)(1-b^2)
        \end{array}
	\right)
$$
The vector $\mathbf{x}_3$ is the rotation axis. 

The expression for $C(-kU)$ after the second equal sign in (\ref{eq:CayM3}) clearly shows that the eigenvectors of 
$U$ and $C(-kU)$ agree in three dimensions as well. 
The general formula for calculating the eigenvalues $\lambda_c$ of $C(-kU)$ from the eigenvalues $\lambda$ of 
$U$ reads:
\begin{equation}
\lambda_{c}=1+\frac{2k\lambda(1+k\lambda)}{1+k^2(a^2+b^2+c^2)}\stackrel{a^2+b^2+c^2=1}{=}
1+\frac{2k\lambda(1+k\lambda)}{1+k^2}.
\label{eq:lambdac3}
\end{equation}
Inserting $\lambda_1$, $\lambda_2$ and $\lambda_3$ in this formula yields:
\begin{equation}
\lambda_{c1}=\cos \vartheta + j \sin \vartheta, \;\;\;
\lambda_{c2}=\cos \vartheta - j \sin \vartheta, \;\;\;
\lambda_{c3}=1.
\end{equation}

We see that in three dimensions the complex eigenvectors (\ref{eq:cevec3}) contain more structure and 
the explicite form of the Cayley transformation (\ref{eq:CayM3}) gets rather unwieldy.

\subsection{Real explanation for three dimensions}

If we follow the treatment of the two-dimensional case given in section \ref{sc:2rex}, then we need to replace the 
imaginary unit $j$ in the eigenvalues $\lambda_1$, $\lambda_2$ and in the eigenvectors $\mathbf{x}_1$, 
$\mathbf{x}_2$ by an element of the real 
three-dimensional geometric algebra $R_3$.\cite{Gu:inotr,DH:CtoG,DH:STC} In principle there are two different choices: The volume element 
$i$ or any two-dimensional unit area element like e.g. 
$\mathbf{i}_1=\mbox{\boldmath $\sigma$}_2\mbox{\boldmath $\sigma$}_3$, 
$\mathbf{i}_2=\mbox{\boldmath $\sigma$}_3\mbox{\boldmath $\sigma$}_1$ or 
$\mathbf{i}_3=\mbox{\boldmath $\sigma$}_1\mbox{\boldmath $\sigma$}_2$. 
($\{\mbox{\boldmath $\sigma$}_1,\mbox{\boldmath $\sigma$}_2,\mbox{\boldmath $\sigma$}_3\}$
denotes an orthonormal basis in $R^3$.)

While both interpretations are possible, let me argue for the second possibility: 
We have seen in section \ref{sc:3comp} that the multiplication of $U$ with 
a vector always projects out the component of this vector parallel to
 $\mathbf{x}_3$ so that the $\mathbf{y}$ on the right hand 
side of equations like $U\mathbf{x} = \mathbf{y}$ is necessarily a vector in the 
two-dimensional plane perpendicular to $\mathbf{x}_3$ 
containing the origin. Thus it seems only natural to interpret the squareroot of 
-1 in the solution of equation (\ref{eq:cpe3}) to be the oriented unit area element 
$
\mathbf{i}=a\mathbf{i}_1+b\mathbf{i}_2+c\mathbf{i}_3
$
 characteristic for the plane 
perpendicular to $\mathbf{x}_3$ containing the origin as opposed to the volume element element $i$ or any other 
two-dimensional unit area element. I will show in the following, that this leads indeed to a consistent interpretation.

Using this area element \textbf{i} we have $\lambda_1=\mathbf{i}$, $\lambda_2=\tilde{\lambda}_1=-\mathbf{i}$ and
\begin{equation}
\mathbf{x}_1=
\left( \begin{array}{c}
        1-a^2 \\
        -ab-\mathbf{i}c \\
        -ac+\mathbf{i}b
       \end{array}\right)
\doteq
\left( \begin{array}{c}
        -ab+\mathbf{i}c \\
        1-b^2 \\
        -bc-\mathbf{i}a
       \end{array}\right)
\doteq
\left( \begin{array}{c}
        -ac-\mathbf{i}b \\
        -bc+\mathbf{i}a \\
        1-c^2
       \end{array}\right), \;\;\;
\mathbf{x}_2= \tilde{\mathbf{x}}_1,
\label{eq:gcevec3}
\end{equation}
where the tilde operation marks the reverse of geometric algebra.
As in the two-dimensional case I interpret the three components of each "eigenvector" 
as spinorial components, i.e. elementary geometric products of two vectors. 
(In the following we will therefore use the expression \textit{three-component eigenspinor} 
instead of "eigenvector".) I again arbitrarily fix one vector from the \textbf{i} plane 
(the plane perpendicular to $\mathbf{x}_3$) as a reference vector \textbf{z} with respect to which I will 
factorize the three component eigenspinors $\mathbf{x}_1$ and $\mathbf{x}_2$. 
With regard to the first representation of the eigenspinors $\mathbf{x}_1$ and $\mathbf{x}_2$ 
we choose to set 
\begin{equation}
\mathbf{z}=\mbox{\boldmath $\sigma$}_{1\|}=\mbox{\boldmath $\sigma$}_1\cdot\mathbf{i}\mathbf{i}^{-1}
=(1-a^2)\mbox{\boldmath $\sigma$}_1-ab\mbox{\boldmath $\sigma$}_2-ac\mbox{\boldmath $\sigma$}_3.
\end{equation}
Using $a^2+b^2+c^2=1$, the square
$
\mathbf{z}^2=1-a^2
$
 is seen to be the first component spinor of the first representation of $\mathbf{x}_1$ and $\mathbf{x}_2$ as given
in (\ref{eq:gcevec3}). 
 
Next we will use the inverse of \textbf{z} 
\begin{equation}
 \mathbf{z}^{-1}=\mbox{\boldmath $\sigma$}_1-\frac{ab}{1-a^2}\mbox{\boldmath $\sigma$}_2
-\frac{ac}{1-a^2}\mbox{\boldmath $\sigma$}_3
\end{equation}
in order to factorize the two other component spinors
$\mathbf{n}_2\mathbf{z}=-ab-\mathbf{i}c$ 
and $\mathbf{n}_3\mathbf{z}=-ac+\mathbf{i}b$
of $\mathbf{x}_1$ in (\ref{eq:gcevec3}) as well. A somewhat cumbersome calculation\footnote{
A good way to speed up and verify such calculations is geometric algebra software, such as~\cite{Cam:MAPLE}
and others.} 
renders
\begin{equation}
\mathbf{n}_2=\mathbf{n}_2\mathbf{z}\mathbf{z}^{-1}=\mbox{\boldmath $\sigma$}_{2\|}, \;\;\;
\mathbf{n}_3=\mathbf{n}_3\mathbf{z}\mathbf{z}^{-1}=\mbox{\boldmath $\sigma$}_{3\|}.
\end{equation}

Summarizing these calculations we have (setting $\mathbf{n}_1=\mathbf{z}=\mbox{\boldmath $\sigma$}_{1\|}$):
\begin{equation}
\mathbf{x}_1=
\left( \begin{array}{c}
        1-a^2 \\
        -ab-\mathbf{i}c \\
        -ac+\mathbf{i}b
       \end{array}\right)
=
\left( \begin{array}{c}
        \mathbf{n}_1\mathbf{z} \\
        \mathbf{n}_2\mathbf{z} \\
        \mathbf{n}_3\mathbf{z}
       \end{array}\right)
=
\left( \begin{array}{c}
        \mbox{\boldmath $\sigma$}_{1\|}\mbox{\boldmath $\sigma$}_{1\|} \\
        \mbox{\boldmath $\sigma$}_{2\|}\mbox{\boldmath $\sigma$}_{1\|} \\
        \mbox{\boldmath $\sigma$}_{3\|}\mbox{\boldmath $\sigma$}_{1\|}
       \end{array}\right).
\end{equation}

The other two equivalent representations of $\mathbf{x}_1$ given in (\ref{eq:gcevec3}) can be written as:
\begin{equation}
\left( \begin{array}{c}
        -ab+\mathbf{i}c \\
        1-b^2 \\
        -bc-\mathbf{i}a
       \end{array}\right)
=
\left( \begin{array}{c}
        \mbox{\boldmath $\sigma$}_{1\|}\mbox{\boldmath $\sigma$}_{2\|} \\
        \mbox{\boldmath $\sigma$}_{2\|}\mbox{\boldmath $\sigma$}_{2\|} \\
        \mbox{\boldmath $\sigma$}_{3\|}\mbox{\boldmath $\sigma$}_{2\|}
       \end{array}\right),
\;\; 
\left( \begin{array}{c}
        -ac-\mathbf{i}b \\
        -bc+\mathbf{i}a \\
        1-c^2
       \end{array}\right)
=
\left( \begin{array}{c}
        \mbox{\boldmath $\sigma$}_{1\|}\mbox{\boldmath $\sigma$}_{3\|} \\
        \mbox{\boldmath $\sigma$}_{2\|}\mbox{\boldmath $\sigma$}_{3\|} \\
        \mbox{\boldmath $\sigma$}_{3\|}\mbox{\boldmath $\sigma$}_{3\|}
       \end{array}\right).
\end{equation}
We see that this simply corresponds to a different choice of the reference vector \textbf{z}, 
as $\mathbf{z'}=\mbox{\boldmath $\sigma$}_{2\|}$ and as $\mathbf{z''}=\mbox{\boldmath $\sigma$}_{3\|}$, 
respectively. In general all possible ways to write $\mathbf{x}_1$ correspond 
to different choices of \textbf{z} from the \textbf{i} plane. 
The geometric product $R_{\mathbf{z}^{-1}\mathbf{z'}} = \mathbf{z}^{-1}\mathbf{z'}$ for any two 
such reference vectors $\mathbf{z}$ and $\mathbf{z'}$ gives the 
rotation operation to rotate one choice of three-component eigenspinor representation $\mathbf{x}_1(\mathbf{z})$ 
 into the other  $\mathbf{x}_1(\mathbf{z'})= \mathbf{x}_1(\mathbf{z})R_{\mathbf{z}^{-1}\mathbf{z'}}$.

As for two dimensions on page \pageref{pg:2rightf} we could also try to interpret $\mathbf{x}_2$ by factoring
out a reference vector \textbf{z} to the right. But since according to equation (\ref{eq:gcevec3}) 
$\mathbf{x}_2$ is simply the reverse of $\mathbf{x}_1$, it seems not really needed for a real interpretation. 
Doing it nevertheless, yields less handy expressions. 

So all we need to give a real geometric interpretation for the three-component eigenspinors $\mathbf{x}_1$ 
(and $\mathbf{x}_2$) is the 
triplet $(\mbox{\boldmath $\sigma$}_{1\|},\mbox{\boldmath $\sigma$}_{2\|},\mbox{\boldmath $\sigma$}_{3\|})$ 
of projections of the three basis vectors $\mbox{\boldmath $\sigma$}_1$, $\mbox{\boldmath $\sigma$}_2$ and 
$\mbox{\boldmath $\sigma$}_3$  onto the \textbf{i} plane. Multiplying this triplet with any vector \textbf{z}, 
element of the \textbf{i} plane,
from the right (from the left) yields all representations of  $\mathbf{x}_1$ (and $\mathbf{x}_2$).

After successfully clarifying the real interpretation of the "complex eigenvectors" in terms of a real vector space 
$R^3$ vector triplet, we turn briefly to the interpretation of the eigenvalues. The real oriented plane unit area
element 
eigenvalues $\lambda_1=\mathbf{i}$ and $\lambda_2=\tilde{\lambda}_1=-\mathbf{i}$ yield via equation 
(\ref{eq:lambdac3}) the eigenvalues of the Cayley transformation $C(-kU)$ as:
\begin{equation}
\lambda_{c1}=\cos \vartheta + \mathbf{i} \sin \vartheta, \;\;\;
\lambda_{c2}=\tilde{\lambda}_{c1}=\cos \vartheta - \mathbf{i} \sin \vartheta, \;\;\;
\lambda_{c3}=1.
\end{equation}

The action of $\lambda_{c1}$ on $\mathbf{x}_1$ and  $\lambda_{c2}$ on $\mathbf{x}_2$, respectively, give
\begin{equation}
C(-kU)\mathbf{x}_1=\lambda_{c1}\mathbf{x}_1=
\lambda_{c1}\left( \begin{array}{c}
        \mathbf{n}_1\mathbf{z} \\
        \mathbf{n}_2\mathbf{z} \\
        \mathbf{n}_3\mathbf{z}
       \end{array}\right)
=
\left( \begin{array}{c}
        \lambda_{c1}\mathbf{n}_1\mathbf{z} \\
        \lambda_{c1}\mathbf{n}_2\mathbf{z} \\
        \lambda_{c1}\mathbf{n}_3\mathbf{z}
       \end{array}\right),
\end{equation}
and
\begin{equation}
C(-kU)\mathbf{x}_2=\lambda_{c2}\mathbf{x}_2=
\lambda_{c2}\left( \begin{array}{c}
        \mathbf{z}\mathbf{n}_1 \\
        \mathbf{z}\mathbf{n}_2 \\
        \mathbf{z}\mathbf{n}_3
       \end{array}\right)
=
\left( \begin{array}{c}
        \mathbf{z}\lambda_{c1}\mathbf{n}_1 \\
        \mathbf{z}\lambda_{c1}\mathbf{n}_2 \\
        \mathbf{z}\lambda_{c1}\mathbf{n}_3
       \end{array}\right).
\label{eq:l2x2int}
\end{equation}
In equation (\ref{eq:l2x2int}) we have used the facts that $\mathbf{x}_2=\tilde{\mathbf{x}}_1$ and that
$\lambda_{c2}\mathbf{z}=\mathbf{z}\tilde{\lambda}_{c2}=\mathbf{z}\lambda_{c1}$, 
since $\mathbf{z}$ is element of the \textbf{i} plane.

We can therefore consistently interpret the both $\lambda_{c1}\mathbf{x}_1$ and $\lambda_{c2}\mathbf{x}_2$
as one and the same rotation of the vector triplet $(\mathbf{n}_1,\mathbf{n}_2,\mathbf{n}_3)
=(\mbox{\boldmath $\sigma$}_{1\|},\mbox{\boldmath $\sigma$}_{2\|},\mbox{\boldmath $\sigma$}_{3\|})$ by the angle
$-\vartheta$ in a right handed sense around the axis $\mathbf{x}_3$ in the \textbf{i} plane. But we are
equally free to alternatively view it as a $+\vartheta$ rotation (in the \textbf{i} plane)
 of the reference vector \textbf{z} instead. No further discussion for the eigenvalues  $\lambda_1$, 
$\lambda_2$ of $U$ is needed, since these are special cases of $\lambda_{c1}$, $\lambda_{c2}$ 
with $\vartheta=\pi/2$.

The third eigenvalue of the Cayley transformation $C(-kU)$ is $\lambda_{c3}=1$, 
which means that any component parallel to $\mathbf{x}_3$ will be invariant under multiplication 
with $C(-kU)$.

So far geometric algebra has served us as an investigative tool in order to gain a consistent real geometric
vector space interpretation of imaginary eigenvalues and complex eigenvectors of antisymmetric matrices 
in two and three dimensions. But it is equally possible to pretend not to know about the antisymmetic matrices
and their eigenvalues and eigenvectors in the first place, and synthetically construct relationships in
geometric algebra which give all the counterparts found in our investigative (analytical) 
appraoach so far. As shown in~\cite{EH:Nagpur} this necessitates in three dimensions the use of the two-sided
spinorial description\cite{DH:NF1,DH:STC} of rotations.

\section{Four Euclidean dimensions - complex treatment}

Even if I don't give the real geometric interpretation for four dimensions, I think it is already worthwhile to
at least define the problem.

Any anti-symmetric matrix in four Euclidean dimensions is proportional to a matrix of the form
\begin{equation}
U=\left(
\begin{array}{cccc}
0 & -e & -f & -g \\
e & 0 & -c & b \\
f & c & 0 & -a \\
g & -b & a & 0
\end{array}
\right)
\end{equation}
with $a^2+b^2+c^2+e^2+f^2+g^2=1$. The characteristic polynomial equation of the matrix $U$ is
\begin{equation}
|U-\lambda E| = 
\left|
\begin{array}{cccc}
-\lambda & -e & -f & -g \\
e & -\lambda & -c & b \\
f & c & -\lambda & -a \\
g & -b & a & -\lambda
\end{array}
\right|
= \lambda^4+\lambda^2+(ae+bf+cg)^2
=0.
\end{equation}
The four eigenvalues are obtained as
\begin{equation}
\lambda_{1,2,3,4}=\pm j \frac{1}{\sqrt{2}} \sqrt{1\pm\sqrt{1-4(ae+bf+cg)^2}}
\end{equation}
One possible representation of the four coordinates of the four complex eigenvectors $\mathbf{x}_n$ 
$(n=1,2,3,4)$ is
\begin{equation}
x_{n1}=\lambda_n(\lambda_n^2+\vec{u}^2), \;\;\;\;\;
\left(
  \begin{array}{c}
     x_{n2} \\
     x_{n3} \\
     x_{n4}
  \end{array}
\right)
=
\lambda_n^2\vec{v}+\lambda_n\;\vec{u}\times\vec{v}+\vec{u}\cdot\vec{v}\;\vec{u}
\end{equation}
where $\vec{u}$ and $\vec{v}$ are three dimensional vectors defined as
$$
\vec{u}=
\left( \begin{array}{c}
        a \\
        b \\
        c
       \end{array}\right),\;\;\;\;
\vec{v}=
\left( \begin{array}{c}
        e \\
        f \\
        g
       \end{array}\right)
$$
and the crossproduct and scalar product are the usual products of three-dimensional 
vector calculus.

The Cayley transformation $C(-kU)$, with the real scalar $k$ is
\begin{eqnarray}
C(-kU) & = & (E+(-kU))^{-1}{}(E-(-kU)) \nonumber \\
& = &  E+2\;\frac{k(1+k^2)U+k^2(1+k^2)U^2+k^3U^3+k^4U^4}{1+k^2+\vec{u}\cdot\vec{v}\;k^4}.
\label{eq:C4}
\end{eqnarray}
The last line of equation (\ref{eq:C4}) shows that $U$ and $C(-kU)$ must have the same
eigenvectors $\mathbf{x}_n$ $(n=1,2,3,4)$. The corresponding eigenvalues of $C(-kU)$ 
can now easily be calculated from (\ref{eq:C4}) as
\begin{equation}
   \lambda_{cn}=1+2\;\frac{k(1+k^2)\lambda_n+k^2(1+k^2)\lambda_n^2+k^3\lambda_n^3+k^4\lambda_n^4}
{1+k^2+\vec{u}\cdot\vec{v}\;k^4}.
\end{equation}

In the future we face the question how to interpret the complex eigenvalues and eigenvectors of
the four-dimensional antisymmetric matrix $U$ and its Cayley transformation $C(-kU)$ in terms of
purely real geometry. I expect to obtain an answer by working with the real geometric algebra $R_4$ of
a real four-dimensional Euclidean vector space $R^4$. This is more involved than for the 
lower dimensional cases. 
The solution will be published elsewhere.

The interpretation of the four dimensional Euclidean problem should also pave the way for treating
the analogous Minkowski space problem. Theoretical physics, especially special relativity, 
electrodynamics and relativistic quantum mechanics may benefit from this.

\subsection*{Acknowledgements}
I first of all thank God for the joy of studying his creation: 
"\ldots since the creation of the world God's invisible qualities - his eternal power and divine nature - 
have been clearly seen, being understood from what has been made \ldots"\cite{P:Rom}. I thank my wife for 
encouragement, T. Ido 
for pointing out my mistakes, H. Ishi for discussions,
 and O. Giering and J.S.R. Chisholm for attracting me 
to geometry. Fukui University provided a good research environment.
I thank K. Shinoda (Kyoto) for his prayerful support.

\end{document}